\documentclass[12pt]{amsart} 
\usepackage{amssymb,amsmath,amscd} 
\theoremstyle{plain} 
\newtheorem{Lem}{Lemma}[section] 
\newtheorem{Prop}[Lem]{Proposition} 
\newtheorem{Thm}[Lem]{Theorem} 
 
\theoremstyle{definition} 
\newtheorem{Def}[Lem]{Definition}

\errorcontextlines=0  \numberwithin{equation}{section}

\newcommand{\calA}{{\mathcal A}} 
\newcommand{\calB}{{\mathcal B}} 
\newcommand{\calH}{{\mathcal H}} 
\newcommand{\calL}{{\mathcal L}} 
\newcommand{\calM}{{\mathcal M}} 
 
\newcommand{\calP}{{\mathcal P}} 
\newcommand{\calS}{{\mathcal S}} 
\newcommand{\calT}{{\mathcal T}} 
\newcommand{\calV}{{\mathcal V}} 
\newcommand{\Id}{\mathrm{Id}}

\newcommand{\fr}{{\mathfrak r}}
\newcommand{\fR}{{\mathfrak R}}
 
\newcommand{\bpi}{\begin{picture}} 
\newcommand{\epi}{\end{picture}}

\newcommand{\regrep}{\breve\omega}

\begin{document} 

\title[Ideals and tensor normal forms] 
{Ideal decompositions and computation of tensor normal forms} 
\date{December 2000} 
\author[B. Fiedler]{Bernd Fiedler}  
\address{Bernd Fiedler\\ Alfred-Rosch-Str. 13\\ 
D-04249 Leipzig \\ Germany}
\urladdr{http://home.t-online.de/home/Bernd.Fiedler.RoschStr.Leipzig/}  
\email{Bernd.Fiedler.RoschStr.Leipzig@t-online.de}  

\begin{abstract}
Symmetry properties of $r$-times covariant tensors $T$ can be described by certain linear subspaces $W$ of the group ring
${\mathbb K} [{\calS}_r]$ of a symmetric group ${\calS}_r$.
If for a class of tensors $T$ such a $W$ is known,
the elements of the orthogonal subspace
$W^{\bot}$ of $W$ within the dual space
${\mathbb K}[{\calS}_r]^{\ast}$ of
${\mathbb K}[{\calS}_r]$
yield linear identities needed for a treatment of the term combination problem for the coordinates of the $T$. We give the structure of these $W$ for every situation which appears in symbolic tensor calculations by computer. Characterizing idempotents of such $W$ can be determined by means of an ideal decomposition algorithm which works in every semisimple ring up to an isomorphism. Furthermore, we use tools such as the Littlewood-Richardson rule, plethysms and discrete Fourier transforms for ${\calS}_r$ to increase the efficience of calculations. All described methods were implemented in a Mathematica package called {\sf PERMS}.
\end{abstract}

\maketitle 

%
%

\section{The Term Combination Problem for Tensors} \label{sect1}%
The use of computer algebra systems for 
symbolic calculations with tensor expressions
is very important in differential geometry,
tensor analysis and 
general relativity theory.
The investigations of this paper\footnote{This paper is a summary of our Habilitationsschrift \cite{fie16}, where all proofs can be found. A part of the proofs were published in earlier papers \cite{fie8,fie14,fie17}, too.
An abridged version of this summary is the paper \cite{fie19}.}
are motivated by the following {\it term combination problem} or {\it normal form problem} which
occurs within such calculations.

Let us consider real or complex linear combinations
\begin{eqnarray}
 \tau & = & \sum_{i = 1}^n {\alpha}_i T_{(i)} \;\;\;,\;\,\;
 {\alpha}_i \in {\mathbb R} , {\mathbb C} \label{equ1}
\end{eqnarray}
of expressions
$T_{(i)}$
which are formed from the coordinates of certain tensors $A$, $B$, $C$,
$\ldots$ by multiplication and, possibly, contractions of some pairs
of indices.
An example of such an expression is
\begin{eqnarray}
 &A_{i a b c}^{}\, A_{\; j k d}^a \, B^{b d}_{\;\;\;e} \, C^{e c} \;\;\;.&
\label{equ2}
\end{eqnarray}
In (\ref{equ2}) we use
Einstein's summation convention.
Further we assume that each of the numbers of $A, B, C, \ldots$ is constant if we run through the set of the $T_{(i)}$.
Now we aim to carry out symbolic
calculations with expressions of the type (\ref{equ1}), (\ref{equ2})
according to the rules of the Ricci calculus.
We assume that
there is a metric tensor
$g$ which allows us to raise or lower indices, for instance
\[
 T_i^{\;j k b} = g^{b c} \, T_{i \; \; \; c}^{\; j k } \; \; \; , \; \; \;
 T_{i a \; \; k }^{\; \; \; j } = g_{a c}^{} \, T_{i \; \; \; k}^{\; c j }
 \; \;.
\]
If now the tensors
$A$, $B$, $C$, $\ldots$ have symmetries and/or fulfil linear identities\footnote{For instance, the Riemann
tensor has the symmetry
$R_{i j k l} = - R_{j i k l} = - R_{i j l k} = R_{k l i j}$ and fulfils
$R_{i j k l} + R_{i k l j} + R_{i l j k} = 0$.
See P. G\"unther \cite{guenth2} for a 
large collection of identities for the curvature tensor, Weyl tensor, etc.
and covariant derivatives of these tensors.
}, then there exist relations between the $T_{(i)}$ in (\ref{equ1}). (We restrict us to linear relations.) Thus the problem arises to detect such relations in sums (\ref{equ1}), generated by symbolic calculations, and to reduce (\ref{equ1}) to 
linear combinations of linearly independent
$T_{(i)}$
({\it normal forms}).

It is well-known that the representation theory of symmetric groups
${\calS}_r$
yields powerful tools to treat this {\it term combination problem}. 
The connection between tensors and the representation theory of
${\calS}_r$
has been considered already in books
 of J.A. Schouten \cite{schout} (1924),
H. Weyl \cite{weyl1} (1939)
and  H. Boerner \cite{boerner} (1955).
In the 1940s Littlewood has developed and used tools such as the 
Littlewood-Richardson rule and plethysms
for the 
investigation of tensors (see S.A. Fulling et al. \cite{full4} (references) and D.E. Littlewood \cite{littlew1} (appendix)).

Applying the same methods, Fulling, King, Wybourne and
Cummins \cite{full4} have calculated large
lists of normal form terms of polynomials of the Riemann curvature tensor and
its covariant derivatives (by means of the program package {\sf Schur} \cite{schur1}).
In their paper \cite{full4} they formulated the following steps to
solve the
above 
term combination problem
for tensors:
\begin{enumerate}
\renewcommand{\labelenumi}{(\alph{enumi})}
\item{Generate the space spanned by the set of homogeneous monomials of a
definite 'order' or 'degree of homogeneity' formed from
the coordinates of tensors of relevance by multiplication and index-pair 
contraction.}
\item{Construct a basis of this space ({\it normal forms}).}
\item{Present an algorithm for expressing an arbitrary element of the space 
in terms of the basis.}
\end{enumerate}
Our present paper yields a method to solve (b) and (c) for
arbitrary tensors.

%
%

\section{Tensors and the Group Ring of a Symmetric Group}
We make use of the following connection
between $r$-times covariant tensors $T \in {\calT}_r V$ over a finite-dimensional ${\mathbb K}$-vector space $V$
and elements of the {\it group ring}
${\mathbb K} [{\calS}_r]$
of a symmetric group ${\calS}_r$ over a field
${\mathbb K} = {\mathbb R}$ or ${\mathbb K} = {\mathbb C}$.\\
\begin{Def}
 Any tensor
 $T \in {\calT}_r V$
 and any $r$-tuple
 $b := (v_1 , \ldots , v_r ) \in V^r$
 of
 $r$
 vectors from
 $V$
 induce a function
 $T_b : {\calS}_r \rightarrow {\mathbb K}$
 according to the rule
 \begin{eqnarray}
T_b (p) & := & T(v_{p(1)} , \ldots , v_{p(r)})\;\;\;,\;\;\;p \in {\calS}_r \,.
\label{eqn136}%
\end{eqnarray}
We identify this function with the group ring element
$T_b := \sum_{p \in {\calS}_r}T_b (p)\,p \in {\mathbb K} [{\calS}_r]$.
\end{Def}

We allow the linear dependence of the $v_i$
and repetitions of vectors in the above $r$-tuple $b$.
Obviously,
two tensors $S , T \in {\calT}_r V$ fulfil $S = T$ iff
$S_b = T_b$ for all
$b \in V^r$.

We try to describe symmetry properties of tensors with the help of the $T_b$ by the following {\it principle}: If a ''class'' of tensors with a certain symmetry property is given, then we search such a linear subspace
$W \subseteq {\mathbb K}[{\calS}_r]$ that contains all $T_b$ of the tensors from the ''class'' being considered. The linear identities which characterize $W$ yield then identities for the coordinates of the $T$ which can be used in the treatment of the term combination problem.

Every
$a = \sum_{p \in {\calS}_r} a(p)\,p \in {\mathbb K} [{\calS}_r]$
acts as so-called {\it symmetry operator}
$a : T \mapsto aT$ on tensors $T \in {\calT}_r V$ if we define
\begin{eqnarray}
(aT)(v_1 , \ldots , v_r) & := &
\sum_{p \in {\calS}_r}\, a(p)\,T(v_{p(1)} , \ldots , v_{p(r)})
\;\;\;,\;\;\; v_i \in V.
\end{eqnarray}
We denote by '$\ast$'
the mapping
$\ast : a = \sum_{p \in {\calS}_r} a(p)\,p \;\mapsto\; a^{\ast} :=
\sum_{p \in {\calS}_r} a(p)\,p^{-1}$.
Furthermore, if $p \in {\calS}_r$ and
$b = (v_1 , \ldots , v_r) \in V^r$,
then we denote by $pb$
the $r$-tuple
$pb :=  (v_{p(1)} , \ldots , v_{p(r)})$.
Many of our calculations are based on\\
\begin{Lem}\hspace{-1mm}\footnote{See B. Fiedler \cite[Sec.III.1]{fie16} and B. Fiedler \cite{fie17}.}
If
$a = \sum_{q \in {\calS}_r} a(q)\,q \in {\mathbb K} [{\calS}_r]$,
$T \in {\calT}_r V$,
$p , q \in {\calS}_r$ and
$b = ( v_1 , \ldots , v_r) \in V^r$,
then we have 
\begin{eqnarray}
T_b (p \circ q) & = & (qT)_b (p) \;=\; T_{pb}(q) 
\label{equ3n2} \\
q(pb) & = & (p \circ q)b  \label{eqn156} \\ %
(aT)_b & = &
T_b \cdot a^{\ast} 
\label{equ3n1} \\
T_b & = & p \cdot T_{pb}\,. \label{equ3n3}
\end{eqnarray}
\end{Lem}
\vspace{0.2cm}

The following symmetry concepts are used for tensors.
(See B. Fiedler \cite[Sec.III.2]{fie16}  and B. Fiedler \cite{fie17}.
See also R. Merris \cite[pp.151,153,157]{merris},
H. Boerner \cite[p.127]{boerner},
G. Eisenreich \cite[p.601]{eisenreich}.)\\
\begin{Def}
\begin{enumerate}
\renewcommand{\labelenumi}{(\alph{enumi})}
\item{
Let $\fr \subseteq {\mathbb K} [{\calS}_r]$ be a right ideal of ${\mathbb K} [{\calS}_r]$ for 
which an $a \in \fr$ and a $T \in {\calT}_r V$ exist such that
$aT \not= 0$. Then the tensor set
${\calT}_{\fr} := \{ aT \;|\; a \in \fr \,,\, T \in {\calT}_r V \}$
is called the
{\it symmetry class}
of tensors defined by $\fr$.
}
\item{
Let $a_1 , \ldots , a_n \in {\mathbb K} [{\calS}_r]$ be a finite set of group ring 
elements. We say that a tensor $T \in {\calT}_r V$ possesses a
{\it symmetry}
defined by $a_1 , \ldots , a_n$ if $\;T$ satisfies
the linear equation system
$\;\;a_i T = 0$, $(i = 1 , \ldots , n)$.
}
\end{enumerate}
\end{Def}
\vspace{0.2cm}

If $e$ is a generating idempotent of a
right ideal
${\mathfrak r} = e \cdot {\mathbb K}[{\calS}_r]$
that defines a symmetry
class ${\calT}_{\fr}$
then ${\calT}_{\fr}$ fulfils
${\calT}_{\fr} = \{ eT \;|\; T \in {\calT}_r V \}$ and
a tensor $T \in {\calT}_r V$ belongs to ${\calT}_{\fr}$ iff
$eT = T$. (See H. Boerner \cite[p.127]{boerner}  or B. Fiedler \cite[Sec.III.2.1]{fie16}.)

Now it can be shown that all $T_b$ of tensors $T$ which have a symmetry (a) or (b) lie in a certain left ideal of ${\mathbb K}[{\calS}_r]$.\\

\begin{Prop}\hspace{-1mm}\footnote{See B. Fiedler \cite{fie17} or
B. Fiedler \cite[Prop. III.2.5, III.3.1, III.3.4]{fie16}.}
\label{satz31}%
Let $e \in {\mathbb K}[{\calS}_r]$ be an idempotent. Then a
$T \in {\calT}_r V$ 
fulfils the condition
$eT = T$
iff
$T_b \;\in\; {\mathfrak l} := {\mathbb K} [{\calS}_r] \cdot e^{\ast}$ for all
$b \in V^r$, i.e.
all $T_b$ of $T$
lie in the left ideal ${\mathfrak l}$ generated by $e^{\ast}$.
\end{Prop}

\begin{Prop}\hspace{-1mm}\footnote{See B. Fiedler \cite{fie17} or
B. Fiedler \cite[Prop. III.3.3, III.3.4]{fie16}.}
\label{satz30}%
Let $a_1 , \ldots , a_m \in {\mathbb K}[{\calS}_r]$ be given group ring 
elements. A
$T \in {\calT}_r V$ satisfies a system of linear
identities
$\;\;a_i T = 0\;,\;(i = 1 , \ldots , m)$,
iff\\
$T_b \;\in\; {\mathfrak l} :=
\{ x \in {\mathbb K} [{\calS}_r] \;|\;
x \cdot a_i^{\ast} = 0 \,,\, i = 1 , \ldots , m \}$
for all $b \in V^r$, i.e.
all $T_b$ of $T$ lie in the
left annihilator ideal ${\mathfrak l}$ of
the set $\{ a_1^{\ast} , \ldots , a_m^{\ast} \}$.
\end{Prop}
\vspace{0.2cm}
The proofs follow easily from (\ref{equ3n1}). A further result is\\
\begin{Prop}\hspace{-1mm}\footnote{See B. Fiedler \cite{fie17} or
B. Fiedler \cite[Prop. III.2.6]{fie16}.}
\label{satz46}%
If $\dim V \ge r$, then every left ideal
${\mathfrak l} \subseteq {\mathbb K}[{\calS}_r]$ fulfils
${\mathfrak l} = {\calL}_{\mathbb K} \{ T_b \;|\;
T \in {\calT}_{{\mathfrak l}^{\ast}} \,,\, b \in V^r \}$.
(Here ${\calL}_{\mathbb K}$ denotes the forming of the linear closure.)
\end{Prop}
\vspace{0.2cm}

If $\dim V < r$, then the $T_b$ of the tensors from
${\calT}_{{\mathfrak l}^{\ast}}$ 
will span only a linear subspace of
${\mathfrak l}$ 
in general.

In the case of tensors $T$ with index contractions the role of the 
$T_b$ is played by certain
sums $\sum_{b \in {\mathfrak B}_{b_0}} {\gamma}_b\,T_b$, which we now define.\\
\begin{Def} \label{defi1}%
Let $g \in {\calT}_2 V$ be a fundamental tensor with arbitrary signature on $V$ and
${\calB} = \{ n_1 , \ldots , n_d \}$ be an orthonormal basis of $V$ with respect
to $g$. Further
let $r$, $l$ be integers with $2 \le 2 l < r$ and
$b_0 = (v_{2 l + 1} , \ldots , v_r ) \in {\calB}^{r - 2 l}$ be a fixed
$(r - 2 l)$-tuple of vectors from ${\calB}$. Then we denote by
${\mathfrak B}_{b_0}$ the set \footnote{We set
$b_0 := \emptyset$ and
${\mathfrak B}_{\emptyset} :=
\bigl\{
(w_1 , w_1 , w_2 , w_2 , \ldots , w_l , w_l)
\in {\calB}^r \;\bigl|\; (w_1 , \ldots , w_l) \in {\calB}^l \bigr\}$
in the case $r = 2 l > 0$.}
of $r$-tuples of basis vectors
\begin{eqnarray*}
{\mathfrak B}_{b_0} & := &
\bigl\{
(w_1 , w_1 , w_2 , w_2 , \ldots , w_l , w_l , v_{2 l + 1} , \ldots , v_r)
\in {\calB}^r \;\bigl|\; (w_1 , \ldots , w_l) \in {\calB}^l \bigr\} \,.
\end{eqnarray*}
Moreover, we set
${\gamma}_b := \prod_{i = 1}^l \,g(w_i , w_i) \;\in \; \{ 1 , -1 \}$
for every $b \in {\mathfrak B}_{b_0}$.
\end{Def}
\vspace{0.2cm}
\begin{Prop}\hspace{-1mm}\footnote{See B. Fiedler \cite{fie17} or
B. Fiedler \cite[Prop. III.3.31]{fie16}.}
\label{satz40}%
Let $T \in {\calT}_r V$ be a tensor and
$g \in {\calT}_2 V$ be a fundamental tensor.
We determine all tensor coordinates with respect to an orthonormal basis
${\calB} = \{ n_1 , \ldots , n_d \}$ of $V$. Let
$b_0 = (n_{i_{2l + 1}} , \ldots , n_{i_r}) \in {\calB}^{r - 2 l}$ be a fixed
$(r - 2 l)$-tuple of basis vectors. Then
\begin{eqnarray}
\sum_{p \in {\calS}_r}\,
(pT)_{j_1 \;\; j_2 \;\;\ldots\;\; j_l \;\; i_{2 l + 1} \ldots i_r}^{
\;\;\; j_1 \;\; j_2 \;\;\ldots\;\; j_l}\,p
& = &
\sum_{b \in {\mathfrak B}_{b_0}} {\gamma}_b\,T_b \,. \label{eqn150}%
\end{eqnarray}
\end{Prop}
\vspace{0.2cm}

Due to Prop. \ref{satz40} the 
$\sum_{b \in {\mathfrak B}_{b_0}} {\gamma}_b\,T_b$ are objects which contain
information about the tensor coordinates of $T$ with $l$ index-pair contractions. In the case of tensors $T$ with index contractions we search for linear subspaces $W \subseteq {\mathbb K}[{\calS}_r]$ which contain all
$\sum_{b \in {\mathfrak B}_{b_0}} {\gamma}_b\,T_b$ for a fixed $b_0$.

The left ideals ${\mathfrak l}$ from the Propositions \ref{satz31},
\ref{satz30} and \ref{satz46} are the simplest examples of linear subspaces $W$ describing tensor symmetries. Before we give such subspaces $W$ in more complicated cases, we will explain how they can be used in the treatment of the term combination problem for tensors.\footnote{See B. Fiedler \cite{fie17} or
B. Fiedler \cite[Sec. III.1, III.4.1]{fie16}.}

\section{The Treatment of the Term Combination Problem}
The term combination problem from Sec. \ref{sect1} can be reformulated as the
problem to find all linear identities between the summands of given tensor 
expressions
\begin{eqnarray}
{\tau}_{i_1 \ldots i_r} & := &
\sum_{p \in {\calP}} c_p\,T_{i_{p(1)} \ldots i_{p(r)}}
\;\;\;,\;\;\; c_p \in {\mathbb K} \,,\, {\calP} \subseteq {\calS}_r
\;\;\;{\rm or} \label{equ203} \\
{\tau}_{i_{2 l + 1} \ldots i_r}
 & := &
\sum_{p \in {\calP}} c_p\,
(pT)_{j_1 \;\; j_2 \;\;\ldots\;\; j_l \;\; i_{2 l + 1} \ldots i_r}^{
\;\;\; j_1 \;\; j_2 \;\;\ldots\;\; j_l}
\;\;\;,\;\;\; c_p \in {\mathbb K} \,,\, {\calP} \subseteq {\calS}_r
\label{equ204}%
\end{eqnarray}
where $T \in {\calT}_{{\mathfrak l}^{\ast}}$ is a tensor from a given symmetry 
class defined by a left ideal ${\mathfrak l}$ (or a right ideal
$\fr = {\mathfrak l}^{\ast}$).
We assume that (\ref{equ203}) and (\ref{equ204}) are results of
symbolic computer calculations. ${\calP}$ is a subset of
permutations, which is determined by the concrete form of the given
expressions (\ref{equ203}) or (\ref{equ204}).

Let
$W \subseteq {\mathfrak l}$
be a linear subspace
which contains all
$T_b$ or $\sum_{b \in {\mathfrak B}_{b_0}} {\gamma}_b\,T_b$ of $T$.
If we consider the {\it orthogonal subspace}
$W^{\bot} := \{ x \in {\mathbb K} [{\calS}_r]^{\star} \;|\;
\forall\,w \in W: \; \langle x , w \rangle = 0 \}$ of
$W$, then 
every $x \in W^{\bot}$ yields a linear identity for the coordinates of $T$ 
since
\begin{eqnarray}
0 & = & \langle x , T_b \rangle \;=\;
\sum_{p \in {\calS}_r} x_p\,T_b(p) \;=\;
\sum_{p \in {\calS}_r} x_p\,T_{i_{p(1)} \ldots i_{p(r)}}\;\;\;\;{\rm or}
\label{eqn166} \\
0 & = & \langle x , \sum_{b \in {\mathfrak B}_{b_0}} {\gamma}_b\,T_b \rangle \;=\;
\underset{p \in {\calS}_r}
{\sum_{{b \in {\mathfrak B}_{b_0}}}}
{\gamma}_b\,T_b(p)\,x_p
 \;=\;
\sum_{p \in {\calS}_r} x_p\,
(pT)_{j_1 \;\;\ldots\;\; j_l \;\; i_{2 l + 1} \ldots i_r}^{
\;\;\; j_1 \;\;\ldots\;\; j_l} \nonumber
\end{eqnarray}
where $x_p := \langle x , p \rangle$, $p \in {\calS}_r$. (The last steps are
correct if all $b$ occurring in (\ref{eqn166}) are $r$-tuples of basis 
vectors of $V$.) Every identity (\ref{eqn166}) can be used to eliminate 
certain summands in (\ref{equ203}), (\ref{equ204}). If $W$ is spanned by all
$T_b$ or $\sum_{b \in {\mathfrak B}_{b_0}} {\gamma}_b\,T_b$ of the tensors 
considered, then $W^{\bot}$ contains all linear identities which are possible
between summands of expressions (\ref{equ203}), (\ref{equ204}) (compare Prop. \ref{satz46}).

If a basis
$\{ h_1 , \ldots , h_k \}$ of $W$ is known, then
the coefficients $x_p$ of the $x \in W^{\bot}$ can be obtained
from the linear equation system
\begin{eqnarray}
\langle x , h_i \rangle \;=\;
\sum_{p \in {\calS}_r} h_i(p)\,x_p & = & 0 \hspace{1cm}(i = 1 , \ldots , k)
\,. \label{eqn168}%
\end{eqnarray}
Thus an important goal is to find such a basis $\{ h_1 , \ldots , h_k \}$
of $W$. (An efficient algorithm for that purpose is given in
Prop. \ref{satz100}.)

Note that (\ref{eqn168}) is a very large system with a
$(k \times r!)$-coefficient matrix, $k = \dim W$. However, since we only need
identities to reduce sums (\ref{equ203}), (\ref{equ204}), we can restrict us 
to solutions of (\ref{eqn168}) which fulfil $x_p = 0$ for
$p \in {\calS}_r \setminus {\calP}$. This reduces the number of columns to
$|{\calP}|$. Furthermore, every of our spaces $W$ is a linear subspace of a
left ideal ${\mathfrak l} = {\mathbb K}[{\calS}_r]\cdot e$. A decomposition
$e = e_1 + \ldots + e_m$ of the generating idempotent into pairwise 
orthogonal, primitive idempotents $e_i$ induces a decomposition
$W = W_1 \oplus\ldots\oplus W_m$ with
$W_i \subseteq {\mathbb K}[{\calS}_r]\cdot e_i$ and a decomposition of the 
tensors $T \in {\calT}_{{\mathfrak l}^{\ast}}$:
$T = e_1^{\ast} T + \ldots + e_m^{\ast} T$. Then we can transform
(\ref{equ203}), (\ref{equ204}) into expressions formed from the
$e_i^{\ast} T$, for instance
\begin{eqnarray}
(\ref{equ204}) \;\;\;\Rightarrow\;\;\; {\tau}_{i_{2 l + 1} \ldots i_r}
 & := &
\sum_{i = 1}^m \sum_{p \in {\calP}} c_p\,
(p(e_i^{\ast}T))_{j_1 \;\; j_2 \;\;\ldots\;\; j_l \;\; i_{2 l + 1} \ldots
i_r}^{\;\;\; j_1 \;\; j_2 \;\;\ldots\;\; j_l} \,, \label{gl1}%
\end{eqnarray}
and use the smaller equation systems of the smaller spaces $W_i$ to determine
linear identities for the coordinates of the $e_i^{\ast} T$. Note, however, that a reduction of expressions such as (\ref{gl1}) by means of identities of the $W_i$ leads to a linear combination of coordinates of the $e_i^{\ast} T$ which cannot be reckoned back into a linear combination of the coordinates of $T$ in general.

\section{The Algorithm for Ideal Decompositions}
Two problems became visible up to now:
We need methods
\begin{enumerate}
\renewcommand{\labelenumi}{(\alph{enumi})}
\item{to determine generating idempotents $e$ for left/right ideals of
${\mathbb K}[{\calS}_r]$ for which such idempotents are unknown (such as in Prop. \ref{satz30}).}
\item{to decompose a given idempotent $e \in {\mathbb K}[{\calS}_r]$ into a sum
$e = e_1 + \ldots + e_m$ of pairwise orthogonal primitive idempotents $e_i$.}
\end{enumerate}
We developed an algorithm which solves these problems. First versions of this algorithm were designed for ${\mathbb K}[{\calS}_r]$ (see B. Fiedler \cite{fie8,fie14}).
But it turned out that this algorithm works even in an arbitrary semisimple ring ${\mathfrak R}$ which fulfils:
\begin{enumerate}
\renewcommand{\labelenumi}{(\Alph{enumi})}
\item{We know explicitly a decomposition
\begin{eqnarray}
\fR \;=\; \bigoplus_{i = 1}^m \, \fR \cdot y_i
& \;\;\;{\rm or}\;\;\; &
\fR \;=\; \bigoplus_{i = 1}^m \, y_i \cdot \fR \label{gl2}%
\end{eqnarray}
of the full ring ${\mathfrak R}$ into minimal left or right ideals generated by known primitive idempotents $y_i$. Pairwise orthogonality of the $y_i$ is not required.}
\item{We know in
${\mathfrak R}$ a method to construct explicitly a solution 
$x \in {\mathfrak R}$ for every equation
\begin{eqnarray}
e \cdot a \cdot x \cdot e \; = \; e & \;\;\;{\rm or}\;\;\; &
e \cdot x \cdot a \cdot e \; = \; e\,, \label{gl4}%
\end{eqnarray}
where $e \in {\mathfrak R}$ is a primitive idempotent and $a \in {\mathfrak R}$ is a ring element with $e \cdot a \not= 0$ or $a \cdot e \not= 0$, respectively.}
\end{enumerate}
We describe now a version (L) of the algorithm for left ideals of a semisimple ring ${\mathfrak R}$. Obviously an analogous version (R) for right ideals can be formulated, too. (See B. Fiedler \cite[Sec.I.2]{fie16} and B. Fiedler \cite{fie8,fie14}.)

A frequent step of the algorithm is the construction of a generating idempotent for a left ideal ${\mathfrak l} = {\mathfrak R} \cdot e \cdot a$, where $e \in {\mathfrak R}$ is a primitive idempotent and $a \in {\mathfrak R}$ is a ring element with $e \cdot a \not= 0$. This is possible by\\
\begin{Prop}\hspace{- 1mm}\footnote{See B. Fiedler \cite[Prop. I.2.1]{fie16}.
Compare B. Fiedler \cite{fie8,fie14}.}
For the above $e , a \in {\mathfrak R}$ there exists such an 
$x \in \fR$ that
\begin{eqnarray}
  e \cdot a \cdot x \cdot e & = & e \;.
\end{eqnarray}
Moreover, the ring element
$e' := x \cdot e \cdot a$
formed with this 
$x$ is an idempotent which generates the minimal left ideal
$\,\fR \cdot e \cdot a$. \label{satz1}%
\end{Prop}
\vspace{0.2cm}

A second construction orthogonalizes idempotents. Let
${\mathfrak l} = {\mathfrak R} \cdot e$ and
$\tilde{{\mathfrak l}} = {\mathfrak R} \cdot \tilde{e}$ be two left ideals with known generating idempotents
$e$ and $\tilde{e}$. Assume that $e$ is primitive and
$e \cdot \tilde{e} \not= e$, i.e. ${\mathfrak l} \not\subseteq \tilde{\mathfrak l}$.
Then the sum ${\mathfrak l} + \tilde{\mathfrak l} = {\mathfrak l} \oplus \tilde{\mathfrak l}$
is direct since the minimality of ${\mathfrak l}$ yields
${\mathfrak l} \cap \tilde{\mathfrak l} = \{0\}$.
Now we search for new generating idempotents
$f$ and $\tilde{f}$ of
${\mathfrak l}$ and $\tilde{\mathfrak l}$ 
which fulfil
$f \cdot \tilde{f} = \tilde{f} \cdot f = 0$.\\
\begin{Thm}\hspace{- 1mm}\footnote{B. Fiedler \cite[Thm. I.2.4]{fie16} 
and B. Fiedler \cite{fie8,fie14}.
The proof uses the fact that the set
$\{ e - x \cdot e + e \cdot x \cdot e \;|\; x \in \fR \}$
is the set of all generating idempotents of ${\mathfrak l} = \fR \cdot e$
(compare D.S. Passman \cite[p.137]{passman2}).
}
The above orthogonalization problem can be solved in 2 steps:
\renewcommand{\labelenumi}{\rm (\roman{enumi})}
 \begin{enumerate}
  \item{We can find such a ring element
   $x \in \fR$ that
   \begin{eqnarray}
    e \cdot (1 - \tilde{e}) \cdot x \cdot e & = & e \;.
   \end{eqnarray}
If we use this $x$ to form
$f := (1 - \tilde{e}) \cdot x \cdot e$,
then $f$ is a generating idempotent of
${\mathfrak l}$ which satisfies
$\tilde{e} \cdot f = 0$.}
  \item{For a result
   $f$
of step
{\rm (i)} there exists an
   $\tilde{x} \in \fR$ such that
   \begin{eqnarray}
    f \cdot (1 - \tilde{e}) \cdot \tilde{x} \cdot f & = & f \;.
   \end{eqnarray}
If we make use of $\tilde{x}$ to form
$\tilde{f} := \tilde{e} - (1 - \tilde{e}) \cdot \tilde{x} \cdot f \cdot \tilde{e}$,
then $\tilde{f}$ is a generating idempotent of
$\tilde{{\mathfrak l}}$ which fulfils
$f \cdot \tilde{f} = \tilde{f} \cdot f = 0$.}
\end{enumerate} \label{theorem5}%
\end{Thm}
Now we can describe our 
{\sc Algorithm (L)}, which allows us to decompose every left ideal
${\mathfrak l} = {\mathfrak R} \cdot a$ with known generating element $a \not= 0$
into a direct sum
${\mathfrak l} = \bigoplus_{i = 1}^m {\mathfrak l}_i$ of minimal left ideals
${\mathfrak l}_i$ explicitly (see B. Fiedler \cite[Sec.I.2]{fie16} and
B. Fiedler \cite{fie8,fie14}).

A multiplication of (\ref{gl2}) by $a$ yields a sum
\begin{eqnarray}
  {\mathfrak l} \;\;\; = \;\;\; \fR \cdot a & = &
\underset{y_i \cdot a \not= 0} { \sum_{{i = 1} }^m}
\,\fR \cdot y_i \cdot a \label{gl3}%
\end{eqnarray}
of minimal left ideals for ${\mathfrak l}$ which however is not direct in general.
Now we can carry out the following steps:
\begin{enumerate}
 \item{The first summand in (\ref{gl3}) is a minimal left ideal. We denote it by ${\mathfrak l}_1$ and we determine a generating idempotent 
$e_1$ of ${\mathfrak l}_1$ by means of Prop. \ref{satz1}.}
\item{We search for the first minimal left ideal
$\fR \cdot y_i \cdot a$ in (\ref{gl3}) which is not contained in
${\mathfrak l}_1$ that means for which
\begin{eqnarray}
 y_i \cdot a \cdot e_1 & \not= & y_i \cdot a \,.
\end{eqnarray}
We denote it by
${\mathfrak l}_2$ 
and we construct a generating idempotent
$e_2$ of ${\mathfrak l}_2$ by means of Prop. \ref{satz1}.
Since
${\mathfrak l}_2$ is minimal and
${\mathfrak l}_2 \not\subseteq {\mathfrak l}_1$ , we obtain
${\mathfrak l}_1 \cap {\mathfrak l}_2 = \{ 0 \}$.
Thus
${\mathfrak l}_1$ and ${\mathfrak l}_2$ form a direct sum
${\tilde{{\mathfrak l}}}_2 := {\mathfrak l}_1 \oplus {\mathfrak l}_2$.
Because
$e_2$ is primitive and
$e_2 \cdot e_1 \not= e_2$,
we can determine new, orthogonal, generating idempotents
${\hat{f}}_1 , f_2$ of ${\mathfrak l}_1 , {\mathfrak l}_2$ 
by means of Theorem \ref{theorem5}.
Then
${\tilde{f}}_2 := {\hat{f}}_1 + f_2$
is a generating idempotent of
${\tilde{{\mathfrak l}}}_2$.}
\item{Now we search for the next minimal left ideal
$\fR \cdot y_i \cdot a$
in (\ref{gl3}) which is
not contained in ${\tilde{{\mathfrak l}}}_2$ that means for which
\begin{eqnarray}
 y_i \cdot a \cdot {\tilde{f}}_2 & \not= & y_i \cdot a \,.
\end{eqnarray}
We denote it by ${\mathfrak l}_3$. We construct a
primitive
generating idempotent
$e_3$ of ${\mathfrak l}_3$ and
pass over to new orthogonal idempotents ${\hat{f}}_2, f_3$ instead of
${\tilde{f}}_2, e_3$.
This leads us to
the left ideal
${\tilde{{\mathfrak l}}}_3 := {\tilde{{\mathfrak l}}}_2 \oplus {\mathfrak l}_3$ which has the generating
idempotent ${\tilde{f}}_3 := {\hat{f}}_2 + f_3$.}
\item{We continue this procedure until we have 
processed 
all left ideals in (\ref{gl3}). The result is a left ideal
${\tilde{{\mathfrak l}}}_n := {{\mathfrak l}}_1 \oplus\ldots\oplus {\mathfrak l}_n$ and a generating idempotent ${\tilde{f}}_n$ of ${\tilde{{\mathfrak l}}}_n$.}
\end{enumerate}
Obviously,
${\tilde{{\mathfrak l}}}_n \subseteq {\mathfrak l}$ since 
every left ideal ${\mathfrak l}_i$ is a summand in \eqref{gl3}. Furthermore, every summand
$\fR \cdot y_i \cdot a$ from \eqref{gl3} which
had been considered
before we had reached
the ideal ${\tilde{{\mathfrak l}}}_n$ is contained in
${\tilde{{\mathfrak l}}}_{n - 1} = {\mathfrak l}_1 \oplus \ldots \oplus {\mathfrak l}_{n - 1} \subseteq 
{\tilde{{\mathfrak l}}}_n$. All other summands $\fR \cdot y_i \cdot a$ of \eqref{gl3}
lie in ${\tilde{{\mathfrak l}}}_n$.
This leads to
${\mathfrak l} \subseteq {\tilde{{\mathfrak l}}}_n$ and
${\mathfrak l} = {\tilde{{\mathfrak l}}}_n$.

According to Theorem \ref{theorem5} (ii), every 
generating idempotent ${\hat{f}}_k$ of ${\tilde{{\mathfrak l}}}_k$
can be
written as ${\hat{f}}_k = (1 - z_k) \cdot {\tilde{f}}_k$
with a $z_k \in \fR$ which we have already determined to carry out
the orthogonalization
$(e_{k+1}, {\tilde{f}}_k) \mapsto (f_{k+1}, {\hat{f}}_k)$. Thus we can write
{\small
 \begin{eqnarray}
  \tilde{f}_n & = & \hat{f}_{n - 1} \,+\, f_n 
  \; = \;
  (1 - z_{n - 1}) \cdot \tilde{f}_{n - 1} \,+\, f_n 
  \; = \;
  (1 - z_{n - 1}) \cdot (\hat{f}_{n - 2} + f_{n - 1}) \,+\,
  f_n \nonumber \\
  & = &
  (1 - z_{n - 1}) \cdot (1 - z_{n - 2}) \cdot \tilde{f}_{n - 2} \,+\,
  (1 - z_{n - 1}) \cdot f_{n - 1} \,+\,
  f_n \nonumber \\
  & \vdots & \nonumber \\
  & = &
  \sum_{k = 1}^{n - 1}
  (1 - z_{n - 1}) \cdot (1 - z_{n - 2}) \cdot \ldots \cdot (1 - z_k) \cdot
  f_k \,+\, f_n \,. \label{equ13}
 \end{eqnarray}
}
 Formula (\ref{equ13}) presents a decomposition
 $\tilde{f}_n = \sum_{k = 1}^n h_k$ of $\tilde{f}_n$
 into summands which fulfil
$h_k := (1 - z_{n - 1}) \cdot (1 - z_{n - 2}) \cdot \ldots
  \cdot (1 - z_k) \cdot f_k \;\in\;  {\mathfrak l}_k$
and
$h_n := f_n \; \in \; {\mathfrak l}_n$.
 Thus,
 $\tilde{f}_n = \sum_{k = 1}^n h_k$
 is the decomposition of
 $\tilde{f}_n$
 which corresponds to
 ${\mathfrak l} = \bigoplus_{k = 1}^n {\mathfrak l}_k$
 and the
 $h_k$
 are pairwise orthogonal generating idempotents of the
 ${\mathfrak l}_k$.

Obviously, the algorithm (L) solves the above problem (b). Furthermore, the algorithm (L) can be extended to left ideals which are non-direct sums
${\mathfrak l} = \sum_{i = 1}^h {\mathfrak R}\cdot a_i$ by applying its steps to the summands of
${\mathfrak l} = \sum_{i = 1}^h \sum_{j = 1}^m {\mathfrak R}\cdot y_j \cdot a_i$.
Likewise, we can construct generating idempotents and decompositions for right ideals ${\mathfrak r} = \sum_{i = 1}^h a_i \cdot {\mathfrak R}$ by the algorithm version (R). If our left/right ideals are intersections
${\mathfrak l} = \bigcap_{i = 1}^h {\mathfrak R}\cdot e_i$ or
${\mathfrak r} = \bigcap_{i = 1}^h e_i \cdot {\mathfrak R}$ of left/right ideals
($e_i$ idempotents), then their right/left annihilator ideals are
${\calA}_r ({\mathfrak l}) =  \sum_{i = 1}^h (1 - e_i) \cdot {\mathfrak R}$
or
${\calA}_l ({\mathfrak r}) =  \sum_{i = 1}^h {\mathfrak R} \cdot (1 - e_i)$, respectively. In this case we can construct a generating idempotent $e$ of 
${\calA}_r ({\mathfrak l})$ or ${\calA}_l ({\mathfrak r})$ by (R) or (L), respectively, and form $e' := 1 - e$ to obtain a generating idempotent $e'$ of
${\mathfrak l}$ or ${\mathfrak r}$. Thus our algorithms solve problem (a) for non-direct sums or intersections of left/right ideals. See B.Fiedler \cite[Sec.I.4]{fie16} or B. Fiedler \cite{fie14} for further details.

\section{Completions of the Decomposition Algorithms}

{\bf The basic assumptions (A) and (B).} Actual decomposition constructions can only be carried out by the algorithms (L) or (R) if our semisimple ring
 ${\mathfrak R}$ fulfils the above assumtions (A) and (B). This is the case for
\begin{enumerate}
\item{the group ring ${\mathfrak R} = {\mathbb K}[{\calS}_r]$ of a symmetric group,}
\item{a ring
${\mathfrak R} = \bigotimes_{i = 1}^m {\mathbb S}_i^{n_i \times n_i}$,
that is an outer direct product of full
$(n_i \times n_i)$-matrix rings over skew fields
${\mathbb S}_i$,}
\item{all semisimple rings
${\mathfrak R}$ for which an isomorphism
$D : {\mathfrak R} \rightarrow {\mathfrak R}' =
\bigotimes_{i = 1}^m {\mathbb S}_i^{n_i \times n_i}$ onto a ring
${\mathfrak R}'$ of the second type is explicitly known.}
\end{enumerate}
(See B. Fiedler \cite[Sec.I.3]{fie16}.)
According to Wedderburn's Theorem every semisimple ring is isomorphic to a ring of type 2. Thus statement 2 means that the decomposition algorithms (L) and (R) work in every semisimple ring up to an isomorphism.

Let us consider
${\mathfrak R} = \bigotimes_{i = 1}^m {\mathbb S}_i^{n_i \times n_i}$.
We denote by $C_{k l} \in {\mathbb S}_i^{n_i \times n_i}$ a matrix in which exactly the element located in the $k$-th row and the $j$-th column is equal to 
$1 \in {\mathbb S}_i$ whereas all other elements vanish.
Then a decomposition of
${\mathfrak R}$
into minimal left/right ideals is given by the decompositions
${\mathbb S}_i^{n_i \times n_i} = \bigoplus_{j = 1}^{n_i} {\mathbb S}_i^{n_i \times n_i} \cdot C_{j j}$
and
${\mathbb S}_i^{n_i \times n_i} = \bigoplus_{j = 1}^{n_i} C_{j j} \cdot {\mathbb S}_i^{n_i \times n_i}$ of the matrix rings into minimal left/right ideals.
Furthermore there exists a
very fast
procedure to solve (\ref{gl4}) in
${\mathfrak R}$.

Since a primitive idempotent $e \in {\mathfrak R}$ has only one non-vanishing block matrix (i.e. $e = (0,\ldots,E,\ldots,0)$, where
$E \in {\mathbb S}_i^{n_i \times n_i}$ is also a primitive idempotent), an equation such as
$e \cdot a \cdot x \cdot e = e$ leads to a single matrix equation
$E \cdot A \cdot X \cdot E = E$. Moreover, $E$ can be written as
$E = f^t \cdot h$ with row vectors $f, h \in {\mathbb S}_i^{n_i}$, where
$h$ is a non-vanishing row of $E$. If we set
$m = h \cdot A$ and determine non-vanishing elements
$m_{j_0}, f_{k_0}$ of $m, f$, then
$X :=  (m_{j_0})^{-1} (f_{k_0})^{-1} C_{j_0 k_0}$ is a solution of
$E \cdot A \cdot X \cdot E = E$, which yields a solution
$x = (0,\ldots,X,\ldots,0)$ of
$e \cdot a \cdot x \cdot e = e$
(B. Fiedler \cite[Sec. I.1.2, I.3.2]{fie16}).
Obviously, this procedure will run very fast on a computer.

For ${\mathbb K}[{\calS}_r]$ the well-known decomposition of
${\mathbb K}[{\calS}_r]$ into minimal left/right ideals by means of {\it Young symmetrizers} guarantees (A). See B. Fiedler \cite{fie8,fie14} for (B).

If for a semisimple ring ${\mathfrak R}$ an above isomorphism
$D: {\mathfrak R} \rightarrow {\mathfrak R}'$ is known (and practicable on a computer),
then every ideal decomposition problem for ${\mathfrak R}$ can be transferred to 
${\mathfrak R}'$ and treated there by the algorithms (L) and (R).

If
${\mathfrak R} = {\mathbb C}[G]$ is the group ring of a finite group $G$, then we have
${\mathbb S}_i = {\mathbb C}$ for all $i$ and
the isomorphism $D$ is called a {\it discrete Fourier transform} for $G$. Explicit algorithms for such Fourier transforms are known at least for {\it abelian groups}, {\it solvable groups}, {\it supersolvable groups} and {\it symmetric groups} (see
M. Clausen und U. Baum \cite{clausbaum1}).\\*[0.3cm]
{\bf Discrete Fourier transforms.}
In group rings
${\mathfrak R} = {\mathbb C}[G]$ 
of large finite groups $G$ (such as $G = {\calS}_r, r \ge 8$), even a single product $a \cdot b$, $\;a , b \in {\mathbb C}[G]$ can lead to high costs in time and computer memory (see B. Fiedler \cite[Sec. I.1.3, I.5.1]{fie16}). Here the use of a discrete Fourier transform
$D : \fR = {\mathbb C} [G] \;\rightarrow\;
\fR' = \bigotimes_{i = 1}^k \, {\mathbb C}^{n_i \times n_i}$ and the transfer of ideal decomposition problems to $\fR'$ is the most important tool to surmount difficulties. Calculations in $\fR'$ have the following advantages:
\begin{enumerate}
\item{Decompositions (\ref{gl2}) and solutions of (\ref{gl4}) can be constructed very fast in $\fR'$.}
\item{
Every product formed during a run of (L) or (R)
contains a factor which is a primitive idempotent $e \in \fR'$.
Since every such $e$ has only 1 non-vanishing block matrix
$E \in {\mathbb C}^{n_i \times n_i}$, i.e.
$e = ( 0 , \ldots , 0 , E , 0 , \ldots , 0 )$,
the costs for every step of (L) or (R) reduce to the costs of calculations in a ring ${\mathbb C}^{n_i \times n_i}$.}
\item{
The algorithms (L) and (R) can be carried out completely within $\fR'$. Only 
input and output data 
have to be mapped between $\fR$ and $\fR'$ by means of $D$
and $D^{-1}$.
Thus, ''less fast'' Fourier transforms can be useful, too.
(See B. Fiedler \cite[Sec.I.5.1]{fie16}.)}
\item{
In $\fR'$ there is a fast construction of bases of linear subspaces $W$,
which we need to form linear equation systems (\ref{eqn168}).}
\end{enumerate}
To describe this construction, we denote by
$C_{i,a} \in  {\mathbb K}^{n \times n}$ that matrix in which the $i$-th row is equal to a given $a \in  {\mathbb K}^n$ whereas all other rows are filled with $0$.\\
\begin{Prop} \label{satz100}%
Let
${\mathfrak l} = {\mathbb K}^{n \times n} \cdot A$
be a minimal left ideal of
${\mathbb K}^{n \times n}$
with known generating element
$0 \not= A \in {\mathbb K}^{n \times n}$
and
$B = [b_{i j}]_{n,n} \not= 0$ be a matrix from
${\mathbb K}^{n \times n}$.
Determine a row $a \not= 0$ of $A$ and 
a parametric form
$\Lambda$ 
of the solution of the linear equation system
\begin{eqnarray}
\sum_{j =1}^n \, b_{i j}\,{\lambda}_j & = & 0 \;\;\;,\;\;\;
i = 1 , \ldots , n 
\;\;\;,\;\;\;{\lambda}_i \in {\mathbb K} \;\;{\rm (unknowns)}\,.
\label{eqn7}
\end{eqnarray}
Then 
${\calB} := 
\left\{ \, B \cdot C_{i,a} \; | \;
\text{$i$ index for which
${\lambda}_i$ is not a parameter in $\Lambda$} \, \right\}$
is a basis of the
${\mathbb K}$-vector space
$B \cdot {\mathfrak l} = B \cdot {\mathbb K}^{n \times n} \cdot A$.
\end{Prop}
\vspace{0.2cm}
See B. Fiedler \cite[Sec.I.1.2]{fie16} for the proof and other fast basis constructions. Spaces with a structure $B \cdot {\mathfrak l}$ are typical examples of spaces $W$ (see Sec. \ref{sect2}, \ref{sect3}). Further, we see that ${\mathfrak l}$ has the basis
${\calB} = \{ C_{i,a} \; | \; i = 1,\ldots,n \}$
if we use the
identity matrix
$B = \Id \in {\mathbb K}^{n \times n}$ for $B$.

For our tensor investigations we need $\fR = {\mathbb K}[{\calS}_r]$.
M. Clausen und U. Baum \cite{clausbaum1,clausbaum2}
developed a very fast Fourier transform for ${\mathbb K}[{\calS}_r]$, which bases on {\it Young's seminormal representation} of ${\calS}_r$ (see also H. Boerner \cite{boerner2} and A. Kerber \cite[Vol.I, p.75,76]{kerber}).
However, since the interpreter {\mbox{\sf Mathematica}} does not allow the full speed and the optimal storage handling of
this ingenious algorithm,
we use {\it Young's natural representation} of
${\calS}_r$ as discrete Fourier transform in our {\mbox{\sf Mathematica}} package {\sf PERMS} \cite{fie10}. (See H. Boerner \cite[pp.102--108]{boerner}, B. Fiedler \cite[Sec.I.5.2]{fie16}.) This implementation works good at least for ${\calS}_r , r \le 8$.\\*[0.3cm]
{\bf Multiplicities.} Obviously, the efficiency of the algorithms (L) and (R) can be improved if we know before a run of (L) or (R) the multiplicities of equivalent minimal left/right ideals
${\mathfrak l}_i$ or ${\mathfrak r}_i$ within decompositions
${\mathfrak l} = \bigoplus_{i = 1}^m {\mathfrak l}_i$ or
${\mathfrak r} = \bigoplus_{i = 1}^m {\mathfrak r}_i$ searched for.
If the algorithms have constructed such a direct sum of minimal left/right ideals of a fixed equivalence class that the number of summands equals the known multiplicity for this class, then the investigation of the remaining ideals of the class can be cancelled. This reduces the calculation time.

In the case of $\fR = {\mathbb K}[{\calS}_r]$
such 
multiplicies can be calculated by means of the {\it irreducible characters} of 
${\calS}_r$, {\it Frobenius reciprocity}, the {\it Littlewood-Richardson rule} and {\it plethysms}.
The determination of the irreducible characters of ${\calS}_r$ is possible by the
{\it Murnaghan-Nakayama formula}. We implemented all these tools in our {\sf Mathematica} package {\sf PERMS} \cite{fie10} (see B. Fiedler \cite[Sec. II.3--II.6]{fie16} for descriptions of implementations). For plethysms we use a very efficient method of F. S\"anger \cite[pp. 29--33]{saenger}.
(See B. Fiedler \cite[Sec.II.6.3]{fie16}.)

\section{Characterizing Left Ideals of Tensor Products} \label{sect2}%
We continue to list linear subspaces
$W \subseteq {\mathbb K}[{\calS}_r]$ describing tensor symmetries.
In the case of tensor products two types of products can be considered:
$T^{(1)} \otimes\ldots\otimes T^{(m)}$ with possibly different $T^{(i)}$ and
$T \otimes\ldots\otimes T$.\\
\begin{Prop}\hspace{- 1mm}\footnote{See B. Fiedler \cite[Sec.III.3.2]{fie16} and B. Fiedler \cite{fie17}.}
Let ${\mathfrak l}_i \subseteq {\mathbb K} [{\calS}_{r_i}]$ $(i = 1,\ldots , m)$ be left ideals
and
$T^{(i)} \in {\calT}_{{\mathfrak l}_i^{\ast}} \subseteq {\calT}_{r_i} V$ be $r_i$-times 
covariant tensors from the symmetry classes characterized by the
${\mathfrak l}_i$.
Consider the product
\begin{eqnarray}
T & := & T^{(1)} \otimes\ldots\otimes T^{(m)} \;\in\; {\calT}_r V
\;\;\;,\;\;\;
r := r_1 + \ldots + r_m \,. \label{eqn120}%
\end{eqnarray}
For every $i$ we define an embedding
\begin{eqnarray}
{\iota}_i : {\calS}_{r_i} \rightarrow {\calS}_r
& \;\;,\;\; &
({\iota}_i s)(k) := 
\left\{
\begin{array}{ll}
{\Delta}_i + s(k - {\Delta}_i) & {\rm if}\;\; r_{i-1} < k \le r_i \\
k & {\rm else}
\end{array}
\right. \label{equ200}%
\end{eqnarray}
where ${\Delta}_i := r_0 + \ldots + r_{i-1}$ and $r_0 := 0$.
Then the $T_b$ of the tensor {\rm (\ref{eqn120})} fulfil
\begin{eqnarray}
 & \;\;\;\;\; &
\forall\, b \in V^r :\;
 T_b \;\in\; {\mathfrak l} \; := \;
{\mathbb K} [{\calS}_r]\cdot{\calL}\bigl\{{\tilde{\mathfrak l}}_1 \cdot\ldots\cdot {\tilde{\mathfrak l}}_m
\bigr\} \; = \;
{\mathbb K} [{\calS}_r]\cdot\bigl({\tilde{\mathfrak l}}_1 \otimes\ldots\otimes
{\tilde{\mathfrak l}}_m \bigr)  \label{eqn122}%
\end{eqnarray}
where ${\tilde{\mathfrak l}}_i := {\iota}_i ({\mathfrak l}_i)$ are the embeddings of
the ${\mathfrak l}_i$ into ${\mathbb K} [{\calS}_r]$ induced by the ${\iota}_i$. If
$\dim V \ge r$, then the above left ideal ${\mathfrak l}$ is
generated by all $T_b \in {\mathbb K} [{\calS}_r]$ which are 
formed from tensor products {\rm (\ref{eqn120})} of arbitrary tensors
$T^{(i)} \in {\calT}_{{\mathfrak l}_i^{\ast}}$.
\end{Prop}
\vspace{0.2cm}

\begin{Prop}\hspace{- 1mm}\footnote{See B. Fiedler \cite[Sec.III.3.2]{fie16} and B. Fiedler \cite{fie17}.}
Let ${\mathfrak l}_0 \subseteq {\mathbb K} [{\calS}_m]$ be a left ideal and
$T \in {\calT}_{{\mathfrak l}_0^{\ast}} \subseteq {\calT}_m V$
be a tensor of order $m$ from the symmetry class
${\calT}_{{\mathfrak l}_0^{\ast}}$.
Consider the product
\begin{eqnarray}
{\hat T} & := & \underbrace{T \otimes\ldots\otimes T}_{n} \;\in\; {\calT}_{m n} V
\,. \label{eqn121}%
\end{eqnarray}
Then all ${\hat T}_b$, $b \in V^{m n}$, lie in the left ideal
\begin{eqnarray}
{\mathfrak l} \; := \;
{\mathbb K} [{\calS}_{m n}]\cdot{\calL}\bigl\{{\mathfrak l}_1 \cdot\ldots\cdot {\mathfrak l}_n \cdot {\mathfrak l}'
\bigr\} \; = \;
{\mathbb K} [{\calS}_{m n}]\cdot\bigl({\mathfrak l}_1 \otimes\ldots\otimes
{\mathfrak l}_n \otimes {\mathfrak l}' \bigr) & & \label{eqn126}%
\end{eqnarray}
where ${\mathfrak l}_i := {\iota}_i ({\mathfrak l}_0)$ are embeddings of
${\mathfrak l}_0$ into ${\mathbb K} [{\calS}_{m n}]$ which are formed by means of 
mappings {\rm (\ref{equ200})} with $r_1 = \ldots = r_n = m$ and
$r = m n$.
Further ${\mathfrak l}'$ denotes the
1-dimensional ideal
${\mathfrak l}' := {\calL}\{ \sum_{q \in Q}\,q \}$
of ${\mathbb K} [Q]$ where $Q \subset {\calS}_{m n}$ is the subgroup
\begin{eqnarray}
Q & := & \Bigl\{
q = {\textstyle
{\binom{k \cdot m - l} {s(k) \cdot m - l}}}
_{\begin{smallmatrix}{1 \le k \le n}\\
{0 \le l \le m - 1}\end{smallmatrix}}
\in {\calS}_{m n} \;\Bigl|\; s \in {\calS}_n \Bigr\}
\;\,\cong\,\; {\calS}_n \,.
\end{eqnarray}
If
$\dim V \ge m \cdot n$, then the above left ideal ${\mathfrak l}$ is
generated by all ${\hat T}_b \in {\mathbb K} [{\calS}_{m n}]$ 
which are formed from tensor products {\rm (\ref{eqn121})} of arbitrary tensors
$T \in {\calT}_{{\mathfrak l}_0^{\ast}}$.
\end{Prop}
\vspace{0.2cm}

Let $\regrep_G : G \rightarrow GL({\mathbb K}[G])$ denote the {\it regular
representation} of a finite group $G$ defined by
$\regrep_g (f) := g \cdot f$, $g \in G$, $f \in {\mathbb K}[G]$. If we use the 
above left ideals ${\mathfrak l}_i$, ${\mathfrak l}_0$, ${\mathfrak l}$
to define subrepresentations
${\alpha}_i := \regrep_{{\calS}_{r_i}} |_{{\mathfrak l}_i}$,
$\alpha := \regrep_{{\calS}_m} |_{{\mathfrak l}_0}$,
$\beta := \regrep_{{\calS}_r} |_{\mathfrak l}$,
then the representation $\beta$ is equivalent to a Littlewood-Richardson 
product or a plethysm\footnote{See the
references \cite{kerber,kerber3,jameskerb,saenger,littlew1,mcdonald,full4} for the Littlewood-Richardson rule and plethysms.},
respectively (see B. Fiedler \cite[Sec.III.3.2]{fie16}):
\begin{eqnarray}
{\mathfrak l} \; {\rm according}\;{\rm to}\;{\rm (\ref{eqn122})}
& \;\;\Longrightarrow\;\; &
\beta \;\sim\;
{\alpha}_1 \,\#\ldots\#\, {\alpha}_m \uparrow {\calS}_r \label{equ201}\\
{\mathfrak l} \; {\rm according}\;{\rm to}\;{\rm (\ref{eqn126})}
& \;\;\Longrightarrow\;\; &
\beta \;\sim\; \alpha \odot [n] \,. \label{equ202}
\end{eqnarray}
These results correspond to statements of S.A. Fulling et al. \cite{full4}.
(\ref{equ201}) and (\ref{equ202}) yield valuable information about 
multiplicities if one wishes to apply
the algorithm 
(L) to ${\mathfrak l}$.
\section{Subspaces Characterizing Tensors with Index Contractions} \label{sect3}%
First we give a universal linear subspace which contains the group ring
elements
$\sum_{b \in {\mathfrak B}_{b_0}} {\gamma}_b\,T_b$ of a tensor $T$ with $l$ 
index contractions for every value of $\dim V$.\\
\begin{Thm} \label{theorem20}%
Let $V, {\calB}, r, l, g, b_0$ have the meaning given in
Def. \ref{defi1} and Prop. 
\ref{satz40}.
Consider the partition
${\lambda}_0 := (2^l\,,\, 1^{r - 2 l}) \vdash r$ and
the lexicographically smallest standard 
tableau $t$ of ${\lambda}_0$. Form the group \footnote{$G$ is a
semidirect product ${\calH}_t \rtimes Q$ and isomorphic to the wreath 
product ${\calS}_2 \wr {\calS}_l$.}
$G := {\calH}_t \cdot Q$ where ${\calH}_t$ is the group of all horizontal 
permutations of $t$ and $Q \subset {\calV}_t$ is the subgroup of all such 
vertical permutations of $t$ which only permute full rows of $t$ with
length $2$. Then every tensor
$T \in {\calT}_{{\mathfrak l}^{\ast}} \subseteq {\calT}_r V$
(${\mathfrak l} = {\mathbb K}[{\calS}_r]\cdot e$, $e$ idempotent) fulfils
\begin{eqnarray}
\sum_{b \in {\mathfrak B}_{b_0}} {\gamma}_b \, T_b & \in &
1_G \cdot {\mathbb K}[{\calS}_r]\cdot e \;\;\;\;,\;\;\;\;
1_G := \sum_{g \in G}\,g \,.
\end{eqnarray}
Furthermore, if $\dim V \ge r - l$, then there is such a
$b_0 \in {\calB}^{r - 2 l}$ that \footnote{Corollary: If $r - 2 l = 0$, then 
the decomposition of $1_G \cdot {\mathbb K}[{\calS}_r]$ into minimal right 
ideals is characterized by a plethysm
$[2] \odot [l] \sim \sum_{\mu \vdash l} [2 \mu ]$. Thus the number $I$
of
linearly independent
{\it invariants} of $T$ is bounded by the sum $M$ of the multiplicities of minimal left 
ideals ${\mathfrak l}_i$ belonging to partitions $2 \mu$, $\mu \vdash l$, in
${\mathbb K}[{\calS}_r]\cdot e = \bigoplus {\mathfrak l}_i$. If $\dim V \ge
l$, then $I = M$. (See B. Fiedler \cite[Sec.III.4.2]{fie16}. Compare S.A. Fulling et al. \cite{full4}.) \label{fussnote}}
\begin{eqnarray}
1_G  \cdot {\mathbb K}[{\calS}_r]\cdot e & = &
{\calL}_{\mathbb K} \Bigl\{ \sum_{b \in {\mathfrak B}_{b_0}} {\gamma}_b \, T_b 
\;\Bigl|\; T \in {\calT}_{{\mathfrak l}^{\ast}} \Bigr\} \,.
\end{eqnarray}
\end{Thm}
The proof can be found in Sec.~III.3.4 of our Habilitationsschrift \cite{fie16}. If $\dim V < r - l$, then the 
$\sum_{b \in {\mathfrak B}_{b_0}} {\gamma}_b \, T_b$ will span only a
linear subspace of $1_G \cdot {\mathbb K} [{\calS}_r] \cdot e$ in general. To 
describe this subspace, we define:\\
\begin{Def} \label{def34}%
If $\lambda = ({\lambda}_1 , \ldots , {\lambda}_k) \vdash r$ is a partition
with length $| \lambda | = k$ and $(v_1 , \ldots , v_k) \in V^k$ is
a $k$-tuple of
vectors, then we denote by
$\langle \lambda ; v_1 , \ldots , v_k \rangle$ or short
$\langle \lambda ; v_i \rangle$
that $r$-tuple from $V^r$ which has the structure
\begin{eqnarray}
\langle \lambda ; v_1 , \ldots , v_k \rangle & := &
(\underbrace{v_1 , \ldots , v_1}_{{\lambda}_1} ,
\underbrace{v_2 , \ldots , v_2}_{{\lambda}_2} , \ldots ,
\underbrace{v_k , \ldots , v_k}_{{\lambda}_k}) \in V^r \,.
\end{eqnarray}
For every 
$b = (v_1 , \ldots , v_r) \in V^r$,
there exists a unique partition
$\lambda \vdash r$ and a permutation $q \in {\calS}_r$ such that $b$ can be 
written as $b = q \langle \lambda ; w_1 , \ldots , w_{| \lambda |} \rangle$
where $w_1 , \ldots , w_{|\lambda |}$ are the pairwise
different, suitably renumbered
vectors from $b$.
We call $\langle \lambda ; w_1 , \ldots , w_{| \lambda |} \rangle$
a {\it grouping}
of $b$ and $\lambda$ the {\it grouping partition}
of $b$, which we also denote by $\lambda = b^{\,\vdash}$.
\end{Def}
\vspace{0.2cm}
\begin{Def} \label{def37}%
\begin{sloppypar}
Let ${\calB}$ be an orthonormal basis with respect to a fundamental tensor 
$g \in {\calT}_2 V$. We call $(n_{i_1} , \ldots , n_{i_{r'}}) \in {\calB}^{r'}$
smaller
than
$(n_{j_1} , \ldots , n_{j_{r'}}) \in {\calB}^{r'}$ if the first
non-vanishing difference $j_k - i_k$ fulfils $j_k - i_k > 0$.
If $\langle \lambda ; w_1 , \ldots , w_{|\lambda |} \rangle$ and
$\langle \lambda ; w'_1 , \ldots , w'_{|\lambda |} \rangle$ are two groupings
of a fixed $r$-tuple $b \in {\calB}^r$ of basis vectors, then we call
$\langle \lambda ; w_1 , \ldots , w_{|\lambda |} \rangle$ smaller
than $\langle \lambda ; w'_1 , \ldots , w'_{|\lambda |} \rangle$
if the $|\lambda |$-tuple
$(n_{i_1} , \ldots , n_{i_{|\lambda |}}) := (w_1 , \ldots , w_{|\lambda |})$
is smaller than the
$|\lambda |$-tuple
$(n_{j_1} , \ldots , n_{j_{|\lambda |}}) := (w'_1 , \ldots , w'_{|\lambda |})$.
\end{sloppypar}
For every $r$-tuple $b \in {\calB}^r$ there exists a
permutation $p_b \in {\calS}_r$
such that $b$ has a representation
$b =  p_b \langle \lambda ; w_1 , \ldots , w_{|\lambda |} \rangle$
where $\langle \lambda ; w_1 , \ldots , w_{|\lambda |} \rangle$ is
the smallest grouping
of $b$ and $\lambda = b^{\,\vdash}$.
We denote by ${\mathfrak p}$
a single-valued mapping
${\mathfrak p}: {\calB}^r \rightarrow {\calS}_r , b \mapsto {\mathfrak p}(b) := p_b$
which assigns exactly one of such permutations $p_b$ to $b$.

Let $b_0 \in {\calB}^{r - 2 l}$ be an $(r - 2 l)$-tuple of vectors from the
basis ${\calB}$.
We denote by ${\Lambda}_{b_0}$ the set
${\Lambda}_{b_0} :=
\{ \lambda \vdash r \;|\;
\exists\,b \in {\mathfrak B}_{b_0} : \; \lambda = b^{\,\vdash} \}$.
Furthermore, we assign to every partition
$\lambda \in {\Lambda}_{b_0}$
the lexicographically smallest standard tableau
$t_{\lambda}$ of $\lambda$
and the set 
${\calM}_{b_0 , \lambda} :=
\{ {\mathfrak p}(b)^{-1} b \in {\calB}^r \;|\;
b \in {\mathfrak B}_{b_0}\;\; {\rm with}\;\; b^{\,\vdash} = \lambda \}$
of such $r$-tuples which are the smallest groupings of the
$b \in {\mathfrak B}_{b_0}$
with grouping partition $\lambda$.
\end{Def}
\vspace{0.2cm}

\begin{Thm} \label{theorem21}%
Let $V, {\calB}, r, l, g, b_0, {\calT}_{{\mathfrak l}^{\ast}}$ and $e$
have the meaning given in Theorem \ref{theorem20}
and ${\mathfrak p}: {\calB}^r \rightarrow {\calS}_r$ be a mapping of the type described 
in Definition \ref{def37}. Then we have
\begin{eqnarray*}
{\calL}_{\mathbb K}\Bigl\{ \sum_{b \in {\mathfrak B}_{b_0}} {\gamma}_b\,T_b
\;\Bigl|\; T \in {\calT}_{{\mathfrak l}^{\ast}} \Bigr\}
 & = &
\sum_{\lambda\in{\Lambda}_{b_0}}\,
\sum_{\langle\lambda ; w_i \rangle\in{\calM}_{b_0 ; \lambda}}
a_{\langle\lambda ; w_i \rangle}\cdot 1_{{\calH}_{t_{\lambda}}}\cdot
{\mathbb K}[{\calS}_r]\cdot e
\end{eqnarray*}
where
\begin{eqnarray*}
a_{\langle \lambda ; w_i \rangle} & := &
\underset{{\mathfrak p}(b)^{-1} b\,=\,\langle \lambda ; w_i
\rangle}
{\sum_{{b \in {\mathfrak B}_{b_0}}}} {\gamma}_b\, {\mathfrak p}(b)^{-1} \,.
\end{eqnarray*}
\end{Thm}
The proof is given in Sec.~III.3.4 of our Habilitationsschrift \cite{fie16}.
Furthermore, we clear in Sec.~III.3.4 to what extent the subspaces from Theorem \ref{theorem21} are independent from the choice of $b_0$.

\section{Concluding Remarks}
In Sec.~III.3.3 of our Habilitationsschrift \cite{fie16} we determined linear subspaces
$W \subseteq {\mathbb K}[{\calS}_r]$ for tensors $T \in {\calT}_r V$ for which the vector space $V$ has a dimension $\dim V < r$. Some details about this case can be found in the paper \cite{fie17}, too.

In Sec.~III.4.2 of the Habilitationsschrift \cite{fie16} we applied Theorem \ref{theorem20} to {\it invariants} (i.e. $r = 2l$ and $b_0 = \emptyset$). Footnote \ref{fussnote} gives one of the results.

We have tested our methods in computer calculations. Among other things, we treated the term combination problem for quadratic monomials in the coordinates of the Riemannian curvature tensor (B. Fiedler \cite[Sec.III.5.1]{fie16}). The correctness of this calculation was controlled by means of the results of S.A. Fulling et al. \cite{full4}. Furthermore, we used Theorem \ref{theorem21} to verify the {\it standard identity}\footnote{See S.A. Amitsur and J. Levitzki \cite{amlev} and S. Bondari \cite{bondari2}.}
$\sum_{p \in {\calS}_{2 n}} \chi (p) A_{p(1)}\cdot\ldots\cdot A_{p(2 n)} = 0$
for $(n \times n)$-matrices $A_i$ in the case $n = 2$ (B. Fiedler \cite[Sec.III.5.2]{fie16}). For this calculation the full algorithm (R) was required.

We implemented our methods in a {\sf Mathematica} package {\sf PERMS}. This package comprises tools for permutation groups and group rings, 
partitions and tableaux,
special idempotents,
characters, the Littlewood-Richardson rule, plethysms,
discrete Fourier transforms and our algorithms (L) and (R). Furthermore, we used the packages {\sf SYMMETRICA} \cite{kerbkohnlas} and {\sf GAP} \cite{Sch95}.



\end{document}